\documentclass[a4paper,10pt,reqno]{article}

\usepackage[utf8]{inputenc}
\usepackage[T1]{fontenc}
\usepackage{lmodern}
\usepackage[english]{babel}
\usepackage{microtype}

\usepackage{amsmath,amssymb,amsfonts,amsthm,esint}
\usepackage{mathtools,accents}
\usepackage{mathrsfs}
\usepackage{aliascnt}
\usepackage{braket}
\usepackage{bm}

\usepackage[a4paper,margin=3cm]{geometry}
\usepackage[citecolor=blue,colorlinks]{hyperref}

\usepackage{enumerate}
\usepackage{xcolor}

\makeatletter
\g@addto@macro\@floatboxreset\centering
\makeatother

\makeatletter
\def\newaliasedtheorem#1[#2]#3{
	\newaliascnt{#1@alt}{#2}
	\newtheorem{#1}[#1@alt]{#3}
	\expandafter\newcommand\csname #1@altname\endcsname{#3}
}
\makeatother

\newtheoremstyle{slanted}{\topsep}{\topsep}{\slshape}{}{\bfseries}{.}{.5em}{}

\theoremstyle{plain}
\newtheorem{theorem}{Theorem}[section]
\newaliasedtheorem{proposition}[theorem]{Proposition}
\newaliasedtheorem{lemma}[theorem]{Lemma}
\newaliasedtheorem{corollary}[theorem]{Corollary}
\newaliasedtheorem{counterexample}[theorem]{Counterexample}

\theoremstyle{definition}
\newaliasedtheorem{definition}[theorem]{Definition}
\newaliasedtheorem{question}[theorem]{Question}
\newaliasedtheorem{openquestion}[theorem]{Open Question}
\newaliasedtheorem{conjecture}[theorem]{Conjecture}

\theoremstyle{remark}
\newaliasedtheorem{remark}[theorem]{Remark}
\newaliasedtheorem{example}[theorem]{Example}

\newcommand{\setR}{\mathbb{R}}

\newcommand{\eps}{\varepsilon}

\let\altphi\phi
\let\phi\varphi
\let\varphi\altphi
\let\altphi\undefined

\newcommand{\abs}[1]{\left\lvert#1\right\rvert}
\newcommand{\norm}[1]{\left\lVert#1\right\rVert}
\newcommand{\seminorm}[1]{\boldsymbol{[} #1 \boldsymbol{]}}

\newcommand{\di}{\mathop{}\!\mathrm{d}}

\newcommand{\loc}{{\rm loc}}

\newcommand{\leb}{\mathscr{L}}

\DeclareMathOperator{\RCD}{RCD}
\newfont{\tmpf}{cmsy10 scaled 2500}

\def\Xint#1{\mathchoice
	{\XXint\displaystyle\textstyle{#1}}%
	{\XXint\textstyle\scriptstyle{#1}}%
	{\XXint\scriptstyle\scriptscriptstyle{#1}}%
	{\XXint\scriptscriptstyle\scriptscriptstyle{#1}}%
	\!\int}
\def\XXint#1#2#3{{\setbox0=\hbox{$#1{#2#3}{\int}$ }
		\vcenter{\hbox{$#2#3$ }}\kern-.6\wd0}}

\def\dashint{\Xint-}

\begin{document}
	
	\title{On the Sobolev space of functions with derivative of logarithmic
		order}
	
	\author{Elia Bru\'e  \thanks{Scuola Normale Superiore, \url{elia.brue@sns.it, }}~~~~~ 
		Quoc-Hung Nguyen \thanks{Scuola Normale Superiore, \url{quochung.nguyen@sns.it, }} 
		\\\\ Scuola Normale Superiore,  Piazza dei Cavalieri 7, I-56100 Pisa, Italy. }
	\maketitle
	
	\begin{abstract}
		Two notions of ``having a derivative of logarithmic order''
		have been studied. They come from the study of regularity of flows and renormalized solutions for the transport and continuity equation associated to weakly differentiable drifts.
	\end{abstract}
	
	\section{Introduction}
	The aim of this note is to study two classes of functions with derivative of logarithmic order in a suitable sense. Such notion comes up naturally in the study of regular Lagrangian flows and renormalized solutions for transport and continuity equation under a Sobolev assumption on the drift (\cite{lions, Ambrosio04, CrippaDeLellis08, Nguyen18-1,BrueNguyen18}). In order to give a concrete idea of this fact we present a formal computation.

     Let us consider a vector field $b:\setR^d\to \setR^d$ in the Euclidean space of dimension $d$ and assume that $b\in W^{1,p}(\setR^d;\setR^d)$ for some $p\geq 1$. Consider the problem
    	\begin{equation}\label{ODE}
     	\begin{dcases}
    	\frac{\di}{\di t} X(t,x)=b(X(t,x)),\\
    	X(0,x)=x, \qquad
    	\forall x\in \mathbb{R}^d.
	    \end{dcases}\tag{ODE}
    	\end{equation}
	The ODE problem in this setting was studied for a first time by Di Perna, Lions \cite{lions} and extended to the BV framework by Ambrosio in \cite{Ambrosio04}. After these two pioneering works this topic has been received a lot attentions becoming a very thriving research field.
	
	Let us now pass to a formal computation. In order to investigate a regularity of the flow map $X_t$ (this question is central in the theory since every regularity of $X$, even very mild, allows compactness theorems, see \cite{BreschJabin15} for a beautiful application to the compressible Navier-Stockes equation) it is natural to take the space derivate in the equation getting
	\begin{equation*}
		\frac{\di}{\di t} \nabla X_t(x)=\nabla b(X_t(x))\cdot \nabla X_t(x).
	\end{equation*}
	Passing to the modulus and integrating we obtain
	\begin{equation*}
	\log\left(|\nabla X_t(x)|\right) \le \int_0^t|\nabla b|(X_s(x))\di s.
	\end{equation*}
	Observe that, under the assumption that $X_t$ is a measure preserving map (that is very natural in this context and crucial to the development of the theory at least in a relaxed form) exploiting the Sobolev regularity of $b$ we conclude that
	\begin{equation}\label{a}
		\log\left(|\nabla X_t(x)|\right)\in L^p(\setR^d),
	\end{equation}
    with a quantitative estimate
    \begin{equation*}
    	\norm{	\log\left(|\nabla X_t|\right)}_{L^p} \le t\norm{\nabla b}_{L^p}.
    \end{equation*}
	Let us point out that our computation is not rigorous since a priory the gradient of $X_t$ does not exist. However from \eqref{a} we learn that the reasonable regularity for a flow map associated to $b\in W^{1,p}(\setR^d;\setR^d)$ is not the classical Sobolev regularity nor a fractional one but something of logarithmic order.
	The first rigorous result in this direction has been proved in \cite{CrippaDeLellis08} where a class of functions is very similar to the one we study in \autoref{section: Nsp} has been introduced.
	
	\paragraph*{}
	
	The paper is organized as follow.
	In \autoref{section Xgammap} we study the class $X^{\gamma,p}$ defined by mean of the Gagliardo-type semi-norm
	\begin{equation}\label{b}
	   \seminorm{f}_{X^{\gamma,p}}:=\left(\int_{B_{1/3}}\int_{\mathbb{R}^d}\frac{|f(x+h)-f(x)|^p}{|h|^d}\frac{1}{\log(1/|h|)^{1-p\gamma}}\di x \di h\right)^{1/p}.
	\end{equation}
     This class was considered for a first time by Leger \cite{LegerFlavien16} in the case $p=2$, (it was defined in a different but equivalent form using the Fourier transformation, see \autoref{subsection: p=2}) with the aim to study regularity and mixing properties of solutions of the continuity equation.
     The authors of the present paper in \cite{BrueNguyen18} considered again the semi-norm \eqref{b} to prove new regularity estimates for the continuity equation proving also their sharpness.
	
	In this first section we prove the Sobolev embedding for the space $X^{\gamma,p}$ in a sharp form, approximation results in the sense of Lusin and the interpolation inequality between spaces $L^p$, $X^{\gamma,p}$ and $W^{s,p}$.
	The techniques we introduce to study such topic are, up to the authors knowledge, new even in the fractional Sobolev setting. We give indeed a new and very simple proof of the fractional Sobolev embedding theorem and of the Lusin's approximation result for function in $W^{s,p}$.
	
	In \autoref{section: Nsp}, we introduced a second class of functions defined à la  Hajlasz (see for instance \cite{Heinonen}). We give a characterization in term of a finiteness of a suitable discrete logarithmic Dirichlet energy (see \eqref{eq: Phi}). We also study weak differentiability properties of functions in this class and we finally establish a link with the first treated space $X^{\gamma,p}$.
	
	\paragraph*{}
	Throughout the present paper we work in the Euclidean space of dimension $d\geq 1$ endowed with the Lebesgue measure $\leb^d$ and the Euclidean norm $|\cdot |$.
	We denote by $B_r(x)$ the ball of radius $r>0$ centered at $x\in\setR^d$. We often write $B_r$ instead of $B_r(0)$. Let us set
	\begin{equation*}
	\dashint_E f\di x=\frac{1}{\leb^d(E)} \int_E f \di x,\qquad \forall~ E\subset \setR^d\  \text{Borel set},
	\end{equation*}
	and
	\begin{equation*}
	Mf(x):=\sup_{r>0} \dashint_{B_r(x)} |f(y)| \di y,
	\qquad
	\forall~x\in \mathbb{R}^d,
	\end{equation*}
	to denote the Hardy-Littlewood maximal function.
	We often use the expression $a\lesssim_c b$ to mean that there exists a universal constant $C$ depending only on $c$ such that $a\leq C b$. The same convention is adopted for $\gtrsim_c$ and $\simeq_c$.

	
	\section{The space $X^{\gamma,p}$}\label{section Xgammap}
	Let us define the first space.
	\begin{definition}\label{def: logSobolev1}
		Let $p\in (0,\infty)$ and $\gamma\in (0,\infty)$ be fixed. We define the semi-norm
		\begin{equation}\label{eq: logSobolev seminorm1}
		\seminorm{f}_{X^{\gamma,p}}:=\left(\int_{B_{1/3}}\int_{\mathbb{R}^d}\frac{|f(x+h)-f(x)|^p}{|h|^d}\frac{1}{\log(1/|h|)^{1-p\gamma}}\di x \di h\right)^{1/p},
		\end{equation}
		and we set
		\begin{equation}\label{eq: space Ap,gamma}
		X^{\gamma,p}:=
		\left\lbrace 
		f\in L^p(\setR^d):\ \seminorm{f}_{X^{\gamma,p}}<\infty
		\right\rbrace
		.
		\end{equation}
	\end{definition}
	It is immediate to verify that $X^{\gamma,p}$ endowed with 
	\begin{equation*}
	\norm{f}_{X^{\gamma,p}}^p:= \norm{f}_{L^p}^p+\seminorm{f}_{X^{\gamma,p}}^p,
	\end{equation*}
	is a Banach space and the semi-norm $\seminorm{\ \cdot\ }_{X^{\gamma,p}}$ is lower semi-continuous with respect to the strong topology of $L^p$.

	Observe that the kernel
	\begin{equation*}
	K_{\gamma}(h):=\mathbf{1}_{B_{1/3}}\frac{1}{|h|^d\log(1/|h|)^{1-p\gamma }},
	\end{equation*}
	appearing in \eqref{eq: logSobolev seminorm1},
	is singular in $\setR^d$ precisely when $\gamma\geq 0$, as a simple computation shows
	\begin{equation*}
	\int_{\setR^d} K_{\gamma}(h)\di h =\int_{\log(3)}^{\infty} t^{p\gamma-1}\di t.
	\end{equation*}
	Therefore, our semi-norm \eqref{eq: logSobolev seminorm1} is not trivial only when $\gamma\geq  0$. Roughly speaking, it has the aim to measure the $L^p$ norm of a derivative of logarithmic order $\gamma$. Indeed, in order to have $\seminorm{f}_{X^{\gamma,p}}<\infty$ it must be
	\begin{equation*}
	\norm{ f(\cdot +h)-f(\cdot)}_{L^p}  \lesssim \frac{1}{\log(1/|h|)^{\gamma}}.
	\end{equation*}
	
	Let us observe the analogies between our space $X^{\gamma,p}$ with the classical Sobolev spaces of fractional order $W^{s,p}$ (see \cite{Adams75} for a reference on this topic).
	Let $p>0$ and $s\in (0,1)$ be fixed, the space $W^{s,p}$ consists of functions $f\in L^p$ such that
	\begin{equation}\label{eq: fractional semi-norm}
	\seminorm{f}_{W^{s,p}}:=
	\left(\int_{\setR^d}\int_{\mathbb{R}^d}\frac{|f(x+h)-f(x)|^p}{|h|^{d+ps}}\di x \di h\right)^{1/p}<\infty,
	\end{equation}
	and it is endowed with the norm
	\begin{equation*}
	\norm{f}_{W^{s,p}}^p:= \norm{f}_{L^p}^p+\seminorm{f}_{W^{p,s}}^p.
	\end{equation*}

	Understanding $X^{\gamma,p}$ as the space of functions in $L^p$ with derivative of logarithmic order $\gamma$ in $L^p$ it is natural to expect the continuous inclusions $X^{\gamma,p} \subset X^{\gamma',p}\subset W^{s,p}$ when $0\leq \gamma \leq \gamma'$ and $s\in (0,1)$. This is the statement of the following proposition whose proof is a simple exercise.
	
	\begin{proposition}\label{prop: simple immersions}
		Let $p\in (0,\infty)$ be fixed. For any $0\leq \gamma\leq \gamma'$ and $s\in (0,1)$ there holds
		\begin{equation*}
		\seminorm{f}_{X^{\gamma,p}}\leq \seminorm{f}_{X^{\gamma',p}},
		\qquad
		\norm{f}_{X^{\gamma,p}}\lesssim_{s,p,\gamma} \norm{f}_{W^{s,p}}.
		\end{equation*}	
	\end{proposition}

	\subsection{The case $p=2$}\label{subsection: p=2}
	In this section we characterize the space $X^{\gamma,2}$ by mean of the Fourier transform. The precise statement is the following.
	
	\begin{theorem}\label{th: fourier characterization}
		Let $\gamma>0$ be fixed.
		For every 
		$f\in L^2(\setR^d)$ it holds
		\begin{equation*}
		\norm{f}_{X^{\gamma,2}}^2\simeq_{d,\gamma} \norm{f}_{L^2}^2+\int_{|\xi|>1} \log(|\xi|)^{2\gamma}|\hat f(\xi)|^2 \di \xi.
		\end{equation*}
	\end{theorem}
	
	\begin{proof} Using the Plancherel's formula we get for any $h\in \mathbb{R}^d$, 
		\begin{equation*}
		\int_{\setR^d} |f(x+h)-f(x)|^2\di x=\int_{\mathbb{R}^d}|e^{ih\cdot \xi}-1|^2 |\hat f(\xi)|^2 \di \xi=2\int_{\mathbb{R}^d}\left(1-\cos (h\cdot \xi)\right) |\hat f(\xi)|^2 \di \xi.
		\end{equation*} 
		Thus, 
		\begin{align}\label{plancherelformual}
		\int_{B_{1/2}}\int_{\mathbb{R}^d}\frac{|f(x+h)-f(x)|^2}{|h|^d\log(1/|h|)^{1-2\gamma}}\di x\di h
		= 2\int_{\mathbb{R}^d}\left[\int_{B_{1/2}}\frac{1-\cos (h\cdot \xi)}{|h|^d\log(1/|h|)^{1-2\gamma}} \di h\right] |\hat f(\xi)|^2 \di \xi.
		\end{align}
		It is enough to show that 
		\begin{equation}\label{z7}
		\int_{B_{1/2}}\frac{1-\cos (h\cdot \xi)}{|h|^d}\frac{1}{\log(1/|h|)^{1-2\gamma}}\di h 	\simeq_{d,\gamma} 
		|\xi|^2,
		\qquad \text{for every}\ |\xi|\leq 10,
		\end{equation}
		and
		\begin{equation}\label{es-ker}
		\int_{B_{1/2}}\frac{1-\cos (h\cdot \xi)}{|h|^d}\frac{1}{\log(1/|h|)^{1-2\gamma}} \di h
		\simeq_{d,\gamma} \log(|\xi|)^{2\gamma},
		\qquad
		\text{for every}\ |\xi|>10.
		\end{equation}
		In order to prove \eqref{z7} we use the elementary inequality $1-\cos(a) \simeq a^2$ for any $a\in (-2,2)$ obtaining
		\begin{equation*}
		\int_{B_{1/2}}\frac{1-\cos (h\cdot \xi)}{|h|^d}\frac{1}{\log(1/|h|)^{1-2\gamma}} \di h 
		\simeq
		\int_{ B_{1/2} }\frac{(h\cdot\xi)^2}{|h|^{d-1}\log(1/|h|)^{1-2\gamma}} \di h
		\simeq_{d,\gamma} |\xi|^2.
		\end{equation*}
		Let us now prove \eqref{es-ker}. First observe that
		\begin{align*}
		\int_{B_{1/2}\setminus B_{1/(10|\xi|)} }\frac{1-\cos (h\cdot \xi)}{|h|^d\log(1/|h|)^{1-2\gamma}} \di h\leq 2 \int_{B_{1/2}\setminus B_{1/(10|\xi|)} }\frac{1}{|h|^d\log(1/|h|)^{1-2\gamma}} \di h\lesssim_d \log(|\xi|)^{2\gamma},
		\end{align*}
		and 
		\begin{align*}
		\int_{B_{1/(10|\xi|)}  }\frac{1-\cos (h\cdot \xi)}{|h|^d\log(1/|h|)^{1-2\gamma}} \di h
		\lesssim & \int_{ B_{1/(10|\xi|)}  }\frac{(h\cdot \xi)^2}{|h|^d\log(1/|h|)^{1-2\gamma}} \di h\\
		\lesssim & |\xi|^2 \int_{ B_{1/(10|\xi|)}  }\frac{1}{|h|^{d-2}\log(10|\xi|)^{1-2\gamma}} \di h\\
		\lesssim_{d,\gamma} & \log(10|\xi|)^{2\gamma-1}
		\lesssim_{d}  \log(|\xi|)^{2\gamma},
		\end{align*}
		where we used again the elementary inequality $1-\cos(a)\lesssim  |a|^2$ for any $a\in\setR$ and the assumption $|\xi|>10$.
		
		Let us now show the converse inequality.
		Using the Coarea formula we can write
		\begin{align*}
		\int_{B_{1/2} }\frac{1-\cos (h\cdot \xi)}{|h|^d\log(1/|h|)^{1-2\gamma}} \di h
		= & \int_{0}^{1/2}\int_{S^{d-1}}\frac{1-\cos (|\xi| r\theta_1)}{r|\log(r)|^{1-2\gamma}} \di \mathcal{H}^{d-1}(\theta)\di r.
		\end{align*}
		It is elementary to see that, there exists a positive constant $C_d$, depending only on $d$, such that for every $\xi\in \setR^d$ and every $r>0$ with  $r|\xi|\ge 1$ it holds
		\begin{align*}
		\int_{S^{d-1}}\left(1-\cos (|\xi| r\theta_1)\right) \di \mathcal{H}^{d-1}(\theta)
		\geq 
		C_d.
		\end{align*}
		Thus, 
		\begin{align*}
		\int_{B_{1/2} }\frac{1-\cos (h\cdot \xi)}{|h|^d}\frac{1}{\log(1/|h|)^{1-2\gamma}} \di h
		\geq & \int_{1/|\xi|}^{1/2}\int_{S^{d-1}}\frac{1-\cos (|\xi| r\theta_1)}{r|\log(r)|^{1-2\gamma}} \di \mathcal{H}^{d-1}(\theta)\di r\\
		\gtrsim_d & \int_{1/|\xi|}^{1/2}\frac{1}{r|\log(r)|^{1-2\gamma}}\di r
		\gtrsim_{d,\gamma} \log(|\xi|)^{1-\gamma}.
		\end{align*}
		Combining the above inequalities, we get \eqref{es-ker}. The proof is complete. 	
	\end{proof}
	
	\subsection{Sobolev and interpolation inequalities}
	In this section we present a new estimate, that could be seen as a log-Sobolev Embedding, introducing a technique that allow us to get a new simple proof of the well-known fraction-Sobolev embedding (see \autoref{prop: fractionalSobolev}). The main result of this section is the following.
	
	\begin{theorem}
		Let $p>0$ and $\gamma>0$ fixed. For any $f\in L^p(\setR^d)$ it holds
		\begin{equation}\label{eq: SobolevEmbedding pointwise}
		|f(x)|^p\log\left(\frac{|f(x)|}{\norm{f}_{L^p(\mathbb{R}^d)}}+2\right)^{p\gamma} \lesssim_{p,\gamma} |f(x)|^p+\int_{B_{1/3}}\frac{|f(x+h)-f(x)|^p}{|h|^d}\frac{1}{\log(1/|h|)^{1-p\gamma}}\di h,
		\end{equation}
		for $\leb^d$-a.e. every $x\in \setR^d$.
		
		In particular the following log-Sobolev inequality holds true for any $f\in X^{\gamma,p}$,
		\begin{equation}\label{eq: SobolevEmbedding}
		\int_{\setR^d} |f(x)|^p\log\left(\frac{|f(x)|}{\norm{f}_{L^p(\mathbb{R}^d)}}+2\right)^{p\gamma} \di x
		\lesssim_{p,\gamma} \norm{f}_{X^{\gamma,p}}^p.
		\end{equation}
	\end{theorem}
	
	\begin{proof} We may assume without loss of generality that $\norm{f}_{L^p}=1$.
		
		Clearly it is enough to show
		\begin{equation}\label{z1}
		|f(x)|^p\log\left(|f(x)|\right)^{p\gamma}\lesssim_{p,\gamma}
		|f(x)|^p+\int_{B_{1/3}}\frac{|f(x+h)-f(x)|^p}{|h|^d}\frac{1}{\log(1/|h|)^{1-p\gamma}}\di h,
		\end{equation}
		for $\leb^d$-a.e. $x\in \set{z\in \mathbb{R}^d:|f(z)|>C_{p,\gamma}}$, where $C_{p,\gamma}>0$ is a fixed constant depending only on $p$ and $\gamma$. 
		
		For any $t\in (0,1/6)$, using the assumption $\norm{f}_{L^p}=1$ we get
		\begin{align*}
		|f(x)|^p & =\dashint_{t^{1/d}<|h|<(2t)^{1/d}} |f(x)|^p \di h\\
		& \lesssim_p \dashint_{t^{1/d}<|h|< (2t)^{1/d}}|f(x+h)|^p\di h+ \dashint_{t^{1/d}<|h|<(2t)^{1/d}}|f(x+h)-f(x)|^p\di h\\
		&\lesssim_p  \frac{1}{t} + \frac{1}{t}\int_{t^{1/d}<|h|<(2t)^{1/d}}|f(x+h)-f(x)|^p\di h.
		\end{align*}
		Here the constant does not depend on dimension $d$ since $\mathcal{L}^d(B_{2^{1/d}})-\mathcal{L}^d(B_{1})\geq \mathcal{L}^d(B_{1})\geq \mathcal{L}^d([0,1]^d)=1. $
		Let us now fix $0<\lambda<1/6$. We integrate the inequality above against a suitable kernel from $\lambda$ to $1/6$ obtaining
		\begin{align*}
		|f(x)|^p\frac{1}{p\gamma}
		&\left(  \log(1/\lambda)^{p\gamma}-\log(6)^{p\gamma}\right)\\
		&=|f(x)|^p\int^{1/6}_{\lambda} \frac{1}{t\log(1/t)^{1-p\gamma}}\di t \\
		&\lesssim_p \int^{1/2}_{\lambda}\frac{1}{t^2\log(1/t)^{1-p\gamma}}\di t
		+\int_{\lambda}^{1/6} \int_{t^{1/d}<|h|<(2t)^{1/d}}\frac{|f(x+h)-f(x)|^p}{t|\log(t)|^{1-p\gamma}} \di h\frac{\di t}{t}\\
		&\lesssim_p \lambda^{-1}\log(1/\lambda)^{p\gamma}+\int_{B_{1/3}}\frac{|f(x+h)-f(x)|^p}{|h|^d}\frac{1}{\log(1/|h|)^{1-p\gamma}}\di h,
		\end{align*}
		so, rearranging the terms, we end up with
		\begin{equation*}
		(|f(x)|^p-C_{p,\gamma}\lambda^{-1})\log(1/\lambda)^{p\gamma}
		\lesssim_{p,\gamma} 
		|f(x)|^p+\int_{B_{1/3}}\frac{|f(x+h)-f(x)|^p}{|h|^d}\frac{1}{\log(1/|h|)^{1-p\gamma}}\di h,
		\end{equation*}
		for any $0<\lambda<1/6$ and for $\leb^d$ a.e. $x\in \setR^d$. Eventually we can choose $\lambda=2C_{p,\gamma}/|f(x)|^p$ when $x\in \set{z\in \mathbb{R}^d:|f(z)|^p>12C_{p,\gamma}}$ and \eqref{z1} immediately follows. The proof is complete.
	\end{proof}
	The just explained strategy could be applied also to the case of fractional Sobolev spaces obtaining the following result.
	
	\begin{proposition}\label{prop: fractionalSobolev}
		Let us fix $s\in (0,1)$ and $p\in (0,d/s)$. We set $p^*:=\frac{dp}{d-sp}$. For any $f\in W^{s,p}(\setR^d)$ the following point-wise inequality holds true
		\begin{equation}\label{z4}
		|f(x)|^{p^*}\lesssim_{p,d,s} \norm{f}_{L^{p^*}}^{p^*-p}\int_{\setR^d} \frac{|f(x)-f(y)|^p}{|x-y|^{d+ps}} \di y,
		\qquad
		\text{for}\ \leb^d\text{-a.e.}\ x\in \setR^d.
		\end{equation}
		In particular we deduce the well-known Sobolev inequality
		\begin{equation}\label{z5}
		\norm{f}_{L^{p^*}}\lesssim_{p,d,s} \norm{f}_{W^{s,p}}.
		\end{equation}	
	\end{proposition}
	
	\begin{proof}
		For any $x\in\setR^d$ and $t>0$ we can write
		\begin{align*}
		|f(x)|^p = & \dashint_{B_{2t}(0)\setminus B_t(0)} |f(x)|^p \di h\\
		\lesssim_p & \dashint_{B_{2t}(0)\setminus B_t(0)} |f(x+h)-f(x)|^p \di h+\dashint_{B_{2t}(x)\setminus B_t(x)} |f(h)|^p \di h.
		\end{align*}
		By mean of H\"older's inequality we get
		\begin{equation*}
		\dashint_{B_{2t}(x)\setminus B_t(x)} |f(h)|^p \di h
		\lesssim_d \norm{f}_{L^{p*}}^p\frac{1}{t^{d-sp}},
		\end{equation*}
		ending up with
		\begin{equation}\label{z3}
		|f(x)|^p\frac{1}{t^{1+sp}}\lesssim_{p,d} \norm{f}_{L^{p*}}^p\frac{1}{t^{d+1}}+\frac{1}{t^{1+sp}}\dashint_{B_{2t}(0)\setminus B_t(0)} |f(x+h)-f(x)|^p \di h.
		\end{equation}
		Let us fix a parameter $\lambda>0$. Integrating \eqref{z3} with respect to $t$ between $\lambda$ and $\infty$ we get
		\begin{equation*}
		|f(x)|^p\frac{1}{\lambda^{sp}}-C_{p,d,s}\norm{f}_{L^{p*}}^p\frac{1}{\lambda^d}\lesssim_{p,d,s} \int_{\setR^d} \frac{|f(x)-f(y)|^p}{|x-y|^{d+ps}} \di y.
		\end{equation*}
		Choosing $\lambda$ such that $\lambda^{d-sp}=2C_{d,p,s}\left( \frac{\norm{f}_{L^{p^*}}}{|f(x)|}\right)^p$ we get \eqref{z4}.
		Integrating \eqref{z4} with respect $x$ over $\setR^d$ we get \eqref{z5}.	
	\end{proof}

	Some remarks are in order.
	\begin{remark}
		The estimate \eqref{eq: SobolevEmbedding pointwise} could be improved in the following way: for any $\gamma>0$, $p>0$ and every $f\in L^{p,\infty}(\setR^d)$ it holds
		\begin{equation}\label{z15}
		|f(x)|^p\log\left(\frac{|f(x)|}{\norm{f}_{L^{p,\infty}}}+2\right)^{p\gamma} \lesssim_{p,\gamma} |f(x)|^p+\int_{B_{1/3}}\frac{|f(x+h)-f(x)|^p}{|h|^d}\frac{1}{\log(1/|h|)^{1-p\gamma}}\di h,
		\end{equation}
		for $\leb^d$-a.e. every $x\in \setR^d$.
		
		Let us explain how to modify the proof of \eqref{eq: SobolevEmbedding pointwise} to get \eqref{z15}. We first assume $p>1$ and we use the inequality
		\begin{equation*}
			\int_E |f(x)| \di x\lesssim_p \leb^d(E)^{1-1/p}\norm{f}_{L^{p,\infty}},\qquad
			\forall E\  \text{Borel subset of}\ \setR^d,
		\end{equation*}
		obtaining
		\begin{align*}
		|f(x)| 
		& \le \dashint_{t^{1/d}<|h|< (2t)^{1/d}}|f(x+h)|\di h+ \dashint_{t^{1/d}<|h|<(2t)^{1/d}}|f(x+h)-f(x)|\di h\\
		&\lesssim_p  \left(\frac{1}{t}\norm{f}_{L^{p,\infty}}\right)^{1/p} + \left(\frac{1}{t}\int_{t^{1/d}<|h|<(2t)^{1/d}}|f(x+h)-f(x)|^p\di h\right)^{1/p}.
		\end{align*}
	   The proof of \eqref{z15} is achieved arguing exactly as in the proof of \eqref{eq: SobolevEmbedding pointwise}.
	   
	   We finally extend \eqref{z15} to every $p\in (0,1]$ by mean of the elementary inequality
	   \begin{align*}
	   |a^\alpha-b^\alpha|\lesssim_{\alpha}  	|a-b|^\alpha,
	   \end{align*}
	   for any $a,b>0$ and $\alpha\in (0,1]$, see \cite[Theorem 2.4]{Daodiaznguyen}.
		\end{remark}
	
	\begin{remark}
		It is worth to remark that \eqref{eq: SobolevEmbedding} and \eqref{eq: SobolevEmbedding pointwise} are dimension free, namely the constant does not depends on $d$. 
	\end{remark}

	Let us now show the sharpness of \eqref{eq: SobolevEmbedding} (and thus of \eqref{eq: SobolevEmbedding pointwise}).
	We prove that
	\begin{equation}\label{z2}
	\seminorm{\mathbf{1}_{B_r(0)}}_{X^{\gamma,p}} \simeq_{p,\gamma} r^d\log(1/r)^{p\gamma},
	\end{equation}
	for any $0<r<1/6$. Thus the function $\mathbf{1}_{B_r}$ saturates \eqref{eq: SobolevEmbedding}. 
	It is enough to show the inequality $\lesssim$ in \eqref{z2}, since the converse is guaranteed by \eqref{eq: SobolevEmbedding}.
	
	Observe that
	\begin{align*}
	\seminorm{\mathbf{1}_{B_r(0)}}_{X^{\gamma,p}}  = & 2\int_{B_r(0)}\int_{\setR^d\setminus B_r(0)} \mathbf{1}_{|x-y|<1/3}\frac{1}{|x-y|^d}\frac{1}{\log(1/|x-y|)^{1-p\gamma}}\di y\di x\\
	= & 2\int_{B_{r}(0)}\int_{B_{1/3}(0)\setminus B_{r}(x)} \frac{1}{|y|^d}\frac{1}{\log(1/|y|)^{1-p\gamma}}\di y\di x.
	\end{align*}
	Take $x\in B_r(0)$, exploiting the inclusion $B(0,r-|x|)\subset B(x,r)$ we deduce
	\begin{align*}
	\int_{B_{1/3}(0)\setminus B_{r}(x)} \frac{1}{|y|^d}\frac{1}{\log(1/|y|)^{1-p\gamma}}\di y\leq &
	\int_{B_{1/3}(0)\setminus B_{r-|x|}(0)} \frac{1}{|y|^d}\frac{1}{\log(1/|y|)^{1-p\gamma}}\di y\\
	 \simeq_{p,\gamma} & \log(1/(r-|x|))^{p\gamma},
	\end{align*}
	thus
	\begin{equation*}
	\seminorm{\mathbf{1}_{B_r(0)}}_{X^{\gamma,p}}  \lesssim_{p,\gamma}	
	\int_{B_r(0)} \log(1/(r-|x|))^{p\gamma} \di x
	\simeq_{d,p,\gamma}  \int_0^r\eps^{d-1} \log(1/\eps)^{p\gamma} \di \eps
	\simeq_{d,p,\gamma}  r^d \log(1/r)^{p\gamma}.
	\end{equation*}
	The function $\mathbf{1}_{B_r(0)}$ is a very natural candidate to show the sharpness of \eqref{eq: SobolevEmbedding}, since in general this class of results are very related to the isoperimetric inequality. 
	However, it is worth to mention that a notion of logarithmic perimeter of order $\gamma>0$ of a set can be obtained writing
	\begin{equation*}
		LP_{\gamma}(E):=\seminorm{\mathbf{1}_E}_{X^{\gamma,1}}.
	\end{equation*}
     We expect that balls are the only minimizers of $LP_{\gamma}$ as happens in the classical and fractional case (see for instance \cite{FuscoMillotMorini11}), but we do not investigate this problem here.

   The last result we present in this section is an interpolation inequality between spaces $L^p(\setR^d)$, $X^{\gamma,p}$ and $W^{s,p}$.
	\begin{proposition} 
		Let $p>0$, $s\in (0,1)$ and $\gamma>0$ be fixed. For any $f\in L^p(\setR^d)$ we have
		\begin{equation}\label{eq: interpolation}
		\seminorm{f}_{X^{\gamma,p}}\lesssim_{p,s,\gamma} \norm{f}_{L^p}\log\left(2+\frac{\norm{f}_{W^{s,p}}}{\norm{f}_{L^p}}\right)^{\gamma}.	
		\end{equation}
    \end{proposition}

\begin{proof}
     Assume without loss of generality $\norm{f}_{L^p}=1$. Let $\lambda\in (0,1/3)$ be fixed, we have
     \begin{align*}
     	\seminorm{f}_{X^{\gamma,p}}^p
     	&=  \int_{ B_{\lambda}(0)}\int_{\setR^d}\frac{|f(x+h)-f(x)|^p}{|h|^d\log(1/|h|)^{1-p\gamma}}\di x\di h
     	+\int_{B_{1/3}(0)\setminus B_{\lambda}(0)}\int_{\setR^d}\frac{|f(x+h)-f(x)|^p}{|h|^d\log(1/|h|)^{1-p\gamma}}\di x\di h \\ &
     	\leq  \frac{\lambda^{ps}}{\log(1/\lambda)^{1-p\gamma}}\seminorm{f}_{W^{s,p}}^p
     	+2^p\norm{f}_{L^p}^p\int_{B_{1/3}(0)\setminus B_{\lambda}(0)}\frac{1}{|h|^d\log(1/|h|)^{1-p\gamma}}\di h\\ &
     	\lesssim_{d,p,\gamma}  \frac{\lambda^{ps}}{\log(1/\lambda)^{1-p\gamma}}\seminorm{f}_{W^{s,p}}^p
     	+\frac{1}{\log(1/\lambda)^{-p\gamma}}
     	\lesssim_{p,\gamma}  \log(1/\lambda)^{p\gamma}(\lambda^{ps}\seminorm{f}_{W^{s,p}}^p+1).
     \end{align*}
     When $\seminorm{f}_{W^{s,p}}\geq 3^s$ we can plug $\lambda=\seminorm{f}_{W^{s,p}}^{-1/s}$ to the previous expression, otherwise we set $\lambda=1/3$, obtaining
     \begin{equation*}
     	\seminorm{f}_{L^p}\lesssim_{p,s,\gamma} \log\left(2+\norm{f}_{W^{s,p}}\right)^{\gamma},
     \end{equation*}
     that is our thesis.
\end{proof}	
	\begin{remark}
		In the particular case $p=2$ the just stated result \eqref{eq: interpolation} could be achieved using the Jensen inequality and the characterization of $X^{\gamma,2}$ by mean of Fourier transform ( see \autoref{th: fourier characterization}).
	\end{remark}

	\subsection{Lusin's estimate}
	It is well-known that a quantitative Lusin's approximation property (see \cite{liu}, \cite{AmbrosioFuscoPallara}) characterizes Sobolev spaces, even in the very abstract setting of measure metric spaces (see for instance \cite{AmbrosioBrueTrevisan17}).
	In this section we study this approximation property for $X^{\gamma,p}$ functions.
	
	Let us introduce the following notation. For any $f\in X^{\gamma,p}$ we define 
	\begin{equation*}
	L_{\gamma,p}f(x):=\left(\int_{B_{1/3}}\frac{|f(x+h)-f(x)|^p}{|h|^d}\frac{1}{\log(1/|h|)^{1-p\gamma}} \di h\right)^{1/p}~~\forall~x\in \mathbb{R}^d, 
	\end{equation*}
	it follows that $L_{\gamma,p}f\in L^p$ with $\norm{f}_{L^p}=\seminorm{f}_{X^{\gamma,p}}$.
	The main result of the section is the following.
	\begin{theorem}\label{Th: Lusin approximation}
		Let $p>0$ and $\gamma>0$ be fixed. For any $f\in X^{\gamma,p}$ there holds
		\begin{equation}\label{eq: LusinLogHolder}
		|f(x)-f(y)|\lesssim_{d,p,\gamma} \log(1/|x-y|)^{-\gamma}\left( L_{\gamma,p}f(x)+L_{\gamma,p}f(y) \right),
		\end{equation}
		for any $x,y\in \setR^d$ such that $|x-y|<\frac{1}{36}$.
	\end{theorem}
	
	\begin{lemma}
		Let $p>0$, $x,y\in \setR^d$ be fixed. For any $f\in L^p$ there holds
		\begin{equation}\label{eq: intermediate}
		|f(x)-f(y)|^p\lesssim_{d,p}
		\dashint_{B_{3r}(0)\setminus B_r(0)} |f(x+h)-f(x)|^p \di h+\dashint_{B_{3r}(0)\setminus B_r(0)} |f(y+h)-f(y)|^p \di h,
		\end{equation}
		for any $r\geq 2|x-y|$.
	\end{lemma}
	
	\begin{proof}
		Let us estimate
		\begin{align*}
		|f(x)-f(y)|^p= &\dashint_{B_{5r/2}(0)\setminus B_{3r/2}(0)} |f(x)-f(y)|^p\di z\\
		&\lesssim_p \dashint_{B_{5r/2}(x)\setminus B_{3r/2}(x)} |f(x)-f(z)|^p\di z
		+\dashint_{B_{5r/2}(x)\setminus B_{3r/2}(x)} |f(z)-f(y)|^p\di z\\
		& \lesssim_d\dashint_{B_{3r}(0)\setminus B_r(0)} |f(x+h)-f(x)|^p \di h+
		\dashint_{B_{5r/2}(x)\setminus B_{3r/2}(x)} |f(z)-f(y)|^p\di z.
		\end{align*}
		Observe that $B_{5r/2}(x)\setminus B_{3r/2}(x)\subset B_{3r}(y)\setminus B_r(y)$ for any $r\geq 2|x-y|$, thus
		\begin{equation*}
		\dashint_{B_{5r/2}(x)\setminus B_{3r/2}(x)} |f(z)-f(y)|^p\di z
		\lesssim_d \dashint_{B_{3r}(0)\setminus B_r(0)} |f(y+h)-f(y)|^p \di h,
		\end{equation*}
		the proof is complete.
	\end{proof}
	
	We are now ready to conclude the proof of \autoref{Th: Lusin approximation}.
	\begin{proof}
		We integrate both sides of \eqref{eq: intermediate} with respect to $r$ against a suitable kernel and we get
		\begin{align*}
		|f(x)-f(y)|^p\int_{2|x-y|}^{1/3} &\frac{1}{r\log(1/r)^{1-p\gamma}} \di r\\ 
		\lesssim_{d,p}&
		\int_{2|x-y|}^{1/3}  \dashint_{B_{3r}(0)\setminus B_r(0)} |f(x+h)-f(x)|^p \di h \frac{\di r}{r \log(1/r)^{1-p\gamma}}\\
		&+
		\int_{2|x-y|}^{1/3}  \dashint_{B_{3r}(0)\setminus B_r(0)} |f(y+h)-f(y)|^p \di h \frac{\di r}{r \log(1/r)^{1-p\gamma}}.
		\end{align*}
		Observe that
		\begin{equation*}
		\int_{2|x-y|}^{1/3}  \dashint_{B_{3r}(0)\setminus B_r(0)} |f(x+h)-f(x)|^p \di h \frac{\di r}{r \log(1/r)^{1-p\gamma}}
		\lesssim_{d,p,\gamma} (L_{\gamma,p}f(x))^p,
		\end{equation*}
		and that
		\begin{align*}
		\int_{2|x-y|}^{1/3} \frac{1}{r\log(1/r)^{1-p\gamma}} \di r &
		=\frac{1}{p\gamma}\left( \log\left(\frac{1}{2|x-y|}\right)^{p\gamma}-\log(3)^{p\gamma}    \right)\\
		&
		\gtrsim_{p\gamma}
		\log\left(\frac{1}{6|x-y|}\right)^{p\gamma}
		\geq 2^{p\gamma}	\log\left(\frac{1}{|x-y|}\right)^{p\gamma},
		\end{align*}
		where in the last passage we have used $|x-y|\leq \frac{1}{36}$.
		The proof is complete.
	\end{proof}
	
	The just described strategy could be used also to prove the standard Lusin's approximation result for $W^{s,p}$ functions with H\"older functions.
	\begin{proposition}
		Let $p\geq 1$ and $s\in (0,1)$ be fixed. For any $f\in  W^{s,p}$ there holds
		\begin{equation*}
		|f(x)-f(y)|\lesssim_{d,s,p}|x-y|^s(D_{s,p}f(x)+D_{s,p}f(y)),
		\qquad
		\forall x,y\in \setR^d,
		\end{equation*}
		where
		\begin{equation*}
		D_{s,p}f(x):=\left(\int_{\setR^d} \frac{|f(x+h)-f(x)|^p}{|h|^{d+ps}}\right)^{1/p}\in L^p(\setR^d).
		\end{equation*}	
	\end{proposition}
	\begin{proof}
		As in the proof of \autoref{Th: Lusin approximation}, we integrate both sides of \eqref{eq: intermediate} with respect to $r$ against a suitable kernel $K(r)$, in this case we should consider $K(r)=\frac{\mathbf{1}_{r>2|x-y|}}{r^{d+ps}}$.
	\end{proof}
	Let us finaly state and prove a simple partial convers implication of \autoref{Th: Lusin approximation}.
	\begin{proposition}	
		Let $p>0$ and $\gamma>0$ be fixed. Let $f\in L^p(\setR^d)$ satisfy
		\begin{equation}\label{eq: lusinLogHolder}
		|f(x)-f(y)|\leq \log(1/|x-y|)^{-\gamma} (g(x)+g(y)),
		\qquad
		\forall x,y\in \setR^d,
		\end{equation}
		for some $g\in L^p(\setR^d)$. Then $f\in X^{\alpha,p}$ for any $\alpha<\gamma$ with estimate
		\begin{equation*}
		\seminorm{f}_{X^{\alpha,p}}^p\lesssim_d \frac{1}{p(\gamma-\alpha)}\norm{g}_{L^p}^p.	
		\end{equation*}
	\end{proposition}
	
	\begin{proof}
		Let us fix $0<\alpha<\gamma$, we estimate
		\begin{align*}
		\seminorm{f}_{X^{\alpha,p}}^p  
		= & \int_{ B_{1/3} } \frac{1}{|h|^d\log(1/|h|)^{1-p\alpha}} \int_0^{\infty} p\lambda^{p-1} \leb^d(\set{x:\ |f(x+h)-f(x)|>\lambda}) \di \lambda\di h\\
		\leq & \int_{ B_{1/3} } \frac{1}{|h|^d\log(1/|h|)^{1-p\alpha}} \int_0^{\infty} p\lambda^{p-1} \leb^d(\set{x:\ g(x)+g(x+h)>\lambda\log(1/|h|)^{\gamma}}) \di \lambda\di h,	
		\end{align*}
		changing variables according to $\lambda\log(1/|h|)^{\gamma}=t$ we get
		\begin{align*}
		\seminorm{f}_{X^{\alpha,p}} 
		\le & \int_{ B_{1/3} } \frac{1}{|h|^d\log(1/|h|)^{1-p(\alpha-\gamma)}} \int_0^{\infty} p t^{p-1} \leb^d(\set{x:\ g(x)+g(x+h)>t}) \di t\di h\\
		\lesssim & \int_{ B_{1/3} } \frac{1}{|h|^d\log(1/|h|)^{1-p(\alpha-\gamma)}}\di h \norm{g}_{L^p}^p
		\simeq_d  \frac{1}{p(\gamma-\alpha)}\norm{g}_{L^p}^p.
		\end{align*}
		The proof is complete.
	\end{proof}
	In order to get a complete characterization of $X^{\gamma,p}$ in term of the Lusin's inequality \eqref{eq: LusinLogHolder} we need to assume something more than the integrability condition $g\in L^p$. We do not porsue here this point.
	
	\section{The space $N^{s,p}$}\label{section: Nsp}
	The aim of this section is to present and study another class of functions with a derivative of logarithmic order. This class comes up naturally in the study of regularity of Lagrangian flows associated to Sobolev vector fields, see \cite{CrippaDeLellis08}.
	\begin{definition}
		Let $p\geq 1$ and $s\in (0,1]$ be fixed. We say that a function  $f\in L^p_{\loc}(\setR^d)$ belongs to $N^{s,p}$ if there exists a positive function $g\in L^p(\setR^d)$ such that
		\begin{equation}\label{eq: Haijalsz log}
		|f(x)-f(y)|\leq
		|x-y|^s\left(\exp\left\lbrace g(x)+g(y)\right\rbrace-1\right),
		\qquad
		\forall
		x,y\in\setR^d.
		\end{equation}
		
		We set $\seminorm{f}_{N^p}:=\inf\{\norm{g}_{L^p}\}$ where the infimum runs over all possible $g$ satisfying \eqref{eq: Haijalsz log}.	
	\end{definition}

     	The $\seminorm{\ \cdot\ }_{N^{1,p}}$ in general is not a semi-norm. It satisfies the triangle inequality and 
     	\begin{equation}\label{eq: almost seminorm property}
     		\seminorm{\lambda f}_{N^{s,p}}\leq |\lambda|\seminorm{f}_{N^{s,p}},
     	\end{equation}
     	as it can be seen using the elementary identity $e^t-1=\sum_{n\geq 1} \frac{t^n}{n!}$, $\forall t\in \setR$. But in general in \eqref{eq: almost seminorm property} the inequality is strict.

	For any $f\in L^p_{\loc}(\setR^d)$ let us define the functional
	\begin{equation}\label{eq: Phi}
	\Phi^*_{s}f(x):=\sup_{r>0}\dashint_{B_r(x)}\log \left( 1+\frac{|f(x)-f(y)|}{r^s}\right) \di y.
	\end{equation}
	Roughly speaking it can be seen as a discrete fractional logarithmic Dirichlet's energy. The aim of the next proposition is to link the condition \eqref{eq: Haijalsz log} to integrability properties of the function $\Phi^*_sf$.
	
	\begin{proposition}\label{prop1}
		Let $p>1$ and $f\in L^p_{\loc}(\setR^d)$ be fixed.
		Then, $f\in N^{s,p}$ if and only if $\Phi^*_{s}f\in L^p(\setR^d)$ and it holds
		\begin{equation}\label{eq: norm equivalence}
		\seminorm{f}_{N^{s,p}}\simeq_{d,p} \norm{\Phi^*_{s}f}_{L^p}.
		\end{equation}	
	\end{proposition}
	\begin{proof}
		For any $f\in N^{s,p}$ it is immediate to see that $\Phi^*_sf\leq 2 Mg$,
		where $M$ is the Hardy-Littlewood maximal function. Thus we get
		\begin{equation}\label{z6}
		\norm{\Phi^*_sf}_{L^p}\lesssim_{d,p} \norm{g}_{L^p},
		\end{equation}
		thanks to the boundedness of $M$ in $L^p$ (see \cite{Stein}).
		In order to achieve the proof of \eqref{eq: norm equivalence} it remain to show the converse of \eqref{z6}.
		
		For any $x,y\in \setR^d$ let us set $r=|x-y|$, we get
		\begin{align*}
		\log\left(1+\frac{|f(x)-f(y)|}{r^s}\right) &=
		\dashint_{B_r(x)} \log\left(1+\frac{|f(x)-f(y)|}{r^s}\right)\di z\\
		& \leq \dashint_{B_r(x)} \log\left(1+\frac{|f(x)-f(z)|}{r^s}\right)\di z
		+ \dashint_{B_r(x)}\log\left(1+\frac{|f(z)-f(y)|}{r^s}\right)\di z\\
		& \leq \Phi^*_sf(x)+\dashint_{B_r(x)}\log\left(1+\frac{|f(z)-f(y)|}{r^s}\right)\di z,
		\end{align*}
		to estimate the last term it is enough to observe that $B(x,r)\subset B(y,2r)$, obtaining
		\begin{equation*}
		\dashint_{B_r(x)}\log\left(1+\frac{|f(z)-f(y)|}{r^s}\right)\di z\lesssim_d \Phi^*_sf(y).
		\end{equation*}
		Thus we end up with
		\begin{equation*}
		\log\left(1+\frac{|f(x)-f(y)|}{|x-y|^s}\right)\lesssim_d \Phi^*_sf(x)+\Phi^*_sf(y),
		\end{equation*}
		that implies $\seminorm{f}_{N^{s,p}}\lesssim_d\norm{\Phi^*_sf(y)}_{L^p}$, and thus \eqref{eq: norm equivalence}. The proof is complete.
	\end{proof}
    Let us point out that the implication $\Phi^*_sf\in L^p(\setR^d)\implies f\in N^{s,p}$ in the case $s=1$ was used in \cite{CrippaDeLellis08} to prove a regularity result for Lagrangian flows.
	Two remarks are in order.
	\begin{remark}
		The assumption $p>1$ in \autoref{prop1} plays a role only in the implication $f\in N^{s,p}\implies \Phi^*_sf\in L^p(\setR^d)$. Indeed in the case $p=1$ only the weaker implication $f\in N^{s,1}\implies \Phi^*_sf\in L^{1,\infty}(\setR^d)$ is available.
	\end{remark}
	\begin{remark}
		For any $f\in L^p(\setR^d)$ and $q\geq 1$ let us consider the functional
		\begin{equation}
		\Phi^*_{s,q}f(x):=\sup_{r>0}\dashint_{B_r(x)}\log \left( 1+\frac{|f(x)-f(y)|}{r^s}\right)^q \di y.
		\end{equation}
		It is immediate to see that
		\begin{equation*}
		\norm{\Phi^*_{s,q}f}_{L^{p/q}}\lesssim_{d,p} \norm{g}_{L^p}^q,
		\qquad
		\text{when}\ 1\leq q<p, 
		\end{equation*}
		and
		\begin{equation*}
		\norm{\Phi^*_{s,p}f}_{L^{1,\infty}}\lesssim_{d,p} \norm{g}_{L^p}^p.
		\end{equation*}
	\end{remark}
	
	Let us recall the definition of weak differentiability.
	\begin{definition}
		We say that a function $f\in L^0(\setR^d)$ is \textit{weakly differentiable} at a point $x\in \setR^d$ if there exists a linear map $L:\setR^d\to \setR$ such that the sequence of functions
		\begin{equation*}
		\frac{f(x+ry)-f(x)-L(ry)}{r|y|}\to 0,
		\end{equation*}
		for $r\to 0$, locally in measure. 
		More precisely
		\begin{equation*}
		\lim_{r\to 0}	\int_{B_R} \abs{\frac{f(x+ry)-f(x)-L(ry)}{r|y|}}\wedge 1 \di y=0,
		\end{equation*}
		for any $0<R<\infty$. We shall denote $L(y):=\nabla f(x)\cdot y$.
	\end{definition}
    It is well-known that a function $f\in L^0(\setR^d)$ is weakly differentiable al $\leb^d$-a.e. $x\in \setR^d$ if and only if it can be approximate with a Lipschitz function in the Lusin's sense. That is to say for any $\eps>0$ there exists a Lipschitz function $g:\setR^d\to \setR$ such that $\leb^d(\set{f\neq g})<\eps$. See \cite{Federer} for a good reference on this topic.
    
    The aim of our next proposition is to study the weakly differentiability property of a function $f\in N^{1,p}$ in a quantitative manner. Precisely we have the following.
	\begin{proposition}\label{prop2}
		 Every $f\in N^{1,p}$ is weakly differentiable at $\leb^d$-a.e. point. Denoting by $\nabla f$ its weak differential we have the following
		\begin{itemize}
			\item[(i)] $\int_{\setR^d} \log\left(1+|\nabla f|\right)^p \di x\lesssim \seminorm{f}_{N^{1,p}}^p$;
			\item[(ii)] for $\leb^d$-a.e. $x\in \setR^d$  there holds
			\begin{equation*}
			\lim_{r\to 0} \dashint_{B_r(0)} \log \left(1+\frac{|f(x+y)-f(x)-\nabla f(x)\cdot y|}{|y|}\right)^p\di y=0.
			\end{equation*}
		\end{itemize}
	\end{proposition}
	\begin{proof}
		It is straightforward to see that $f$ is weakly differentiable (see discussion below).
		For any constant $M>0$ we have
		\begin{equation*}
		\dashint_{B_1(0)}\log \left( 1+\frac{|f(x+ry)-f(x)|}{r|y|}\wedge M \right) \di y
		\lesssim \Phi^*_1f(x),
		\end{equation*}
		recalling that $y\to \frac{f(x+ry)-f(x)}{r|y|}$ converges locally in measure to $\nabla f(x)\cdot y/|y|$ for $\leb^d$-a.e. $y\in \setR^d$ we deduce
		\begin{equation}\label{z8}
		\dashint_{B_1(0)}\log(1+|\nabla f(x)\cdot (y/|y|)|\wedge M)\di y
		=\lim_{r\to 0} \dashint_{B_1(0)}\log \left( 1+\frac{|f(x+ry)-f(x)|}{r|y|}\wedge M \right) \di y
		\lesssim \Phi^*_1f(x),
		\end{equation}
		for $\leb^d$-a.e. $x\in \setR^d$. It is immediate to deduce (i) from \eqref{z8} and \autoref{prop1}. Let us pass to the proof of (ii).		
	  
	 First of all let us consider an increasing convex function $\Psi$ such that
	  \begin{equation*}
	  	\lim_{t\to \infty}\frac{\Psi(t)}{t}=\infty,
	  \end{equation*}
	  and 
	  \begin{equation*}
	  	\int_{\setR^d} \Psi(g(x)^p)\di x<\infty,
	  \end{equation*}
	  it exists thanks to the well-known Dunford-Pettis lemma, see \cite{AmbrosioFuscoPallara}. We can also assume that $t\mapsto	\frac{\Psi(t)}{t}$ is increasing.
	  Setting
	  \begin{equation*}
	  	f_r(y)=\frac{|f(x+ry)-f(x)-\nabla f(x)\cdot ry|}{r|y|},
	  \end{equation*}
	  and using the very definition of $N^{1,p}$ we get
	  \begin{align*}
	  	\sup_{r>0}\int_{B_1(0)} \Psi(\log(1+f_r(y))^p) \di y\lesssim_{\Psi} & \int_{B_1(0)}\Psi((g(x+ry)+g(x))^p) \di y+\Psi(\log(1+|\nabla f(x)|)^p)\\
	  	\lesssim_{\Psi} & M \left(\Psi(g^p)\right)(x)+\Psi(\log(1+|\nabla f(x)|)^p),
	  \end{align*}
	  thus it is finite for $\leb^d$-a.e. $x\in \setR^d$ thanks to the (1,1) weak estimate for the maximal function.
	  thus it is enough to prove that (ii) holds for every $x\in \setR^d$ such that
	  \begin{equation}
	  	\sup_{r>0}\int_{B_1(0)} \Psi(\log(1+f_r(y))^p) \di y=: T<\infty.
	  \end{equation}
	  
	  Let us fix a parameter $0<\lambda<1/2$, using the Jensen inequality we get
	  \begin{align*}
	  &	\dashint_{B_r(0)} \log \left(1+\frac{|f(x+y)-f(x)-\nabla f(x)\cdot y|}{|y|}\right)^p\di y\\
	  	& ~~~~=\dashint_{B_1(0)} \log(1+f_r(y))^p \di y\\
	  & ~~~\simeq_d \int_{B_1\cap \set{f_r>\lambda^{-1}}} \log(1+f_r(y))^p \di y+ \int_{B_1\cap \set{\lambda<f_r<\lambda^{-1}}} \log(1+f_r(y))^p \di y\\&~~~~+
	  	    \int_{B_1\cap \set{f_r\leq \lambda}} \log(1+f_r(y))^p \di y\\ & ~~~\lesssim_d \frac{\log(1+\lambda^{-1})^p}{\Psi\left(\log(1+\lambda^{-1})^p\right)} T+ \log(1+\lambda^{-1})^p\leb^d(B_1\cap\{f_r>\lambda\})+\log(1+\lambda)^p.
	  \end{align*}
	  Taking first the limit for $r\to 0$ and after $\lambda\to 0$ we get the sought conclusion since $f_r\to 0$ locally in measure.
	\end{proof}

\begin{remark}
	 It is natural to wonder if the statement of \autoref{prop2} has a converse, or if some quantitative version of the weakly differentiability at $\leb^d$-a.e. almost every point could guaranteed the property \eqref{eq: Haijalsz log}. The answer is negative, indeed for any $p\in [1,\infty)$ it is possible to built a function $f\in L^p(\setR^d)$ that is weakly differentiable almost everywhere with $\nabla f=0$ but does not belong in $N^{1,q}$ for any $q$.
	 
	 Let us illustrate how to built such example. Let us fix an integer $M>0$. It is enough to built a function $f$ supported on the set $[0,1]\subset \setR^d$ that is weakly differentiable with $f'=0$ $\leb^1$-a.e., $\norm{f}_{L^{\infty}}=1$ and $\seminorm{f}_{N^{1,1}}\geq cM$, where $c$ does not depend on $f$.
	 Let us define 
	 \begin{equation*}
	 	 f(x)=\sum_{k=0}^{M-1}(-1)^k\mathbf{1}_{(K/M,(K+1)/M]}(x)~~\forall~~x\in [0,1]. 
	 \end{equation*}
	It is trivial to see that $\norm{f}_{L^{\infty}}=1$ and $f$ is differentiable at every point out side a finite set with derivative equal to zero. Let us show that $\seminorm{f}_{N^{1,1}}\geq cM$. Take any $g$ realizing the identity in \eqref{eq: Haijalsz log}, we have
	 \begin{equation}\label{z10}
	 	g(x)\mathbf{1}_{(i/M, (i+1)/M]}(x)+	g(y)\mathbf{1}_{(j/M, (j+1)/M]}(y)\geq \log(1+2/|x-y|),
	 \end{equation}
	 for any $i$ even and $j$ odd. Summing up in \eqref{z10} we get the inequality
	 \begin{equation*}
	 	g(x)+g(y)\geq \frac{M}{2} \log(1+2/|x-y|),
	 \end{equation*}
	 that integrated over $[0,1]\times [0,1]$ trivially implies our thesis.
\end{remark}

    The last result of this section concerns with a link between the spaces $X^{\gamma,p}$ and $N^{s,p}$. A version of this result was crucial in our previous paper \cite{BrueNguyen18} to study a sharp regularity for the continuity equation associated to a divergence-free Sobolev drift.
	\begin{theorem} 
		Let $p\geq 1$ and $s\in (0,1]$ be fixed. For any $f$ satisfying \eqref{eq: Haijalsz log} the following estimate holds true
		\begin{equation}\label{es1}
		\int_{B_{1/3}} \int_{\setR^d} \frac{1\wedge|f(x+h)-f(x)|^q}{|h|^d}\frac{1}{\log(1/|h|)^{1-p}} \di x \di h
		\lesssim_{s,p,d} \norm{g}_{L^p}^p+\norm{g}_{L^q}^q.
		\end{equation}
		In particular we have the continuous immersion
		\begin{equation*}
		N^{s,p}\cap L^{\infty}(\setR^d)\hookrightarrow X^{1,p},	
		\end{equation*}
		for any $s\in (0,1]$ and $p\geq 1$.
	\end{theorem}
	
	\begin{proof}
		By \eqref{eq: Haijalsz log}, we can write 
		\begin{align*}
		&\int_{B_{1/e}} \int_{\setR^d} \frac{1\wedge|f(x+h)-f(x)|^q}{|h|^d\log(1/|h|)^{1-p}} \di x \di h\\
		& \leq \int_{B_{1/e}} \int_{0}^1q\lambda^{q-1}\leb^{d}(\{x:|h|^s\left(\exp\left\lbrace g(x+h)+g(x) \right\rbrace-1\right)>\lambda\}) \di \lambda \frac{1}{|h|^d\log(1/|h|)^{1-p}} \di h\\
		& = \int_{B_{1/e}} \int_{0}^{1/|h|^s}q\lambda^{q-1}\leb^{d}(\{x:\left(\exp\left\lbrace g(x+h)+g(x) \right\rbrace-1\right)>\lambda\}) \di \lambda \frac{1}{|h|^{d-sq}\log(1/|h|)^{1-p}} \di h.
		\end{align*}
		Note that for $0<\lambda<2$, one has 
		\begin{align*}
		  \leb^{d}(\{x:\exp\left\lbrace g(x+h)+g(x) \right\rbrace-1>\lambda\})
		  &\leq \leb^{d}(\{x:4(g(x+h)+g(x) )>\lambda\})\\ 
		  &\leq  2\leb^{d}(\{x:8g(x)>\lambda\}).
		\end{align*}
		Thus, 
		\begin{align*}
		 &\int_{B_{1/e}} \int_{\setR^d} \frac{1\wedge|f(x+h)-f(x)|^q}{|h|^d\log(1/|h|)^{1-p}} \di x \di h\\
		 &\leq \int_{B_{1/e}} \int_{0}^{2}2q\lambda^{q-1}\leb^{d}(\{x:8g(x)>\lambda\})\di \lambda \frac{1}{|h|^{d-sq}\log(1/|h|)^{1-p}} \di h\\
		 & +\int_{B_{1/e}} \int_{2}^{1/|h|^s}q\lambda^{q-1}\leb^{d}(\{x:\exp\left\lbrace g(x+h)+g(x) \right\rbrace >\lambda\}) \di \lambda \frac{1}{|h|^{d-sq}\log(1/|h|)^{1-p}} \di h\\
	    &\lesssim_{p,d,s,q} \int_{0}^{1}q\lambda^{q-1}\leb^d(\{g>\lambda\})\di \lambda
	    +\int_{B_{1/e}} \int_{2}^{1/|h|^s}\lambda^{q-1}\leb^{d}(\{g >\frac{1}{2}\log(\lambda)\}) \di \lambda \frac{1}{|h|^{d-sq}\log(1/|h|)^{1-p}} \di h\\
		& \lesssim_{p,d,s,q}  \norm{g}_{L^q}^q
		+\int_{2}^{\infty}\left[\int_{|h|^s<1/\lambda} \frac{1}{|h|^{d-sq}\log(1/|h|)^{1-p}} \di h\right]\lambda^{q-1} \leb^q(\{g>\frac{1}{2}\log(\lambda)\}) \di \lambda.
		\end{align*}
		Since for $\lambda>2$, 
		\begin{align*}
		\int_{|h|^s<1/\lambda} \frac{1}{|h|^{d-sq}\log(1/|h|)^{1-p}} \di h\lesssim_{p,d,s} \lambda^{-q}\log(\lambda)^{p-1},
		\end{align*} 
		we deduce
		\begin{align*}
		&\int_{B_{1/e}} \int_{\setR^d} \frac{1\wedge|f(x+h)-f(x)|}{|h|^d}\frac{1}{\log(1/|h|)^{1-p}} \di x \di h
		\\
		&\lesssim_{p,d,s,q}\norm{g}_{L^q}^q
		+\int_{1}^{\infty}\lambda^{-1}\log(\lambda)^{p-1} \leb^d(\{ g>\frac{1}{2}\log(\lambda)\}) \di \lambda\\
		&\lesssim_{p,d,s,q} \norm{g}_{L^q}^q
		+\int_{1}^{\infty}\lambda^{p-1} \leb^d(\{g >\lambda\}) \di \lambda,
		\end{align*}
		which implies \eqref{es1}. The proof is complete. 
	\end{proof}
     Let us finally remark that \eqref{es1} could not be improved in the  following way
	\begin{equation}\label{z13}
	  \int_{B_{1/3}} \int_{\setR^d} \frac{|f(x+h)-f(x)|^q}{|h|^d\log(1/|h|)^{1-p}} \di x \di h
	  \lesssim_{p,q,d} (\norm{g}_{L^p}^p+\norm{g}_{L^q}^q)\norm{f}_{L^{\alpha}}^{\beta},
	\end{equation}
	for some $\alpha, \beta>0$.
	Indeed, using \eqref{eq: almost seminorm property} and applying a scaling argument we get $\beta<\beta+\min\{p,q\}\leq q$. Moreover, exploiting the elementary inequality
	\begin{equation*}
		\lambda(e^a-1)\le e^{\lambda^{\eps}C_{\eps}a}-1,
	\end{equation*}
	for any $a\geq 0$ and $\lambda>1$ we obtain that $\seminorm{\lambda f}_{N^{s,q}}^q\leq C_{\eps}\lambda^{q\eps} \seminorm{f}_{X^{s,q}}^q$ for any $s\in (0,1]$ when $\lambda>1$. Plugging $\lambda f$ and the just mentioned estimate in $\eqref{z13}$ we deduce $\beta \geq q$.

\end{document}